\documentclass[twocolumn,aps,pre,floatfix,superscriptaddress,longbibliography]{revtex4-1}
\pdfoutput=1 %This is for the ArXiv submission only
\usepackage{amsmath,amssymb,eucal,graphicx,float,epstopdf}

\usepackage[colorlinks=true, urlcolor=blue, anchorcolor=blue, citecolor=blue,filecolor=blue,linkcolor=blue,menucolor=blue]{hyperref}

\begin{document}
\title{Random Sequential Covering}

\author{P. L. Krapivsky}
\affiliation{Department of Physics, Boston University, Boston, Massachusetts 02215, USA}
\affiliation{Santa Fe Institute, Santa Fe, New Mexico 87501, USA}

\begin{abstract} 
In random sequential covering, identical objects are deposited randomly, irreversibly, and sequentially; only attempts  increasing the coverage are accepted. A finite system eventually gets congested, and we study the statistics of congested configurations. For the covering of an interval by dimers, we determine the average number of deposited dimers, compute all higher cumulants, and establish the probabilities of reaching minimally and maximally congested configurations. We also investigate random covering by segments with $\ell$ sites and sticks. Covering an infinite substrate continues indefinitely, and we analyze the dynamics of random sequential covering of $\mathbb{Z}$ and $\mathbb{R}^d$. 
\end{abstract}

\maketitle

\section{Introduction}
\label{sec:Intro}

The covering of Euclidean spaces by spheres, or other identical objects, is an old and fascinating subject \cite{VDW56,Rogers,GL87,Conway}. A natural guess is that the least dense sphere covering can be obtained by properly enlarging the radii of the densest sphere packing is true \cite{Kershner39,Verblunsky} in two dimensions, but false in three and eight dimensions \cite{Vallentin05,Vallentin06}, and probably false in all $d\geq 3$ with the possible exception of $d=24$. The least dense sphere coverings are more complicated than the densest sphere packings; the latter are also poorly understood and currently known \cite{Hales05,Vallentin16,Cohn17} only when $d=1,2,3,8, 24$. The least dense sphere coverings initially attracted interest thanks to their beauty and intriguing connections with number theory \cite{Rogers,GL87,Conway,Janson86,Kahane,Kendall13,Calka,Penrose21}. There are also intriguing applications to genomics \cite{Athreya04}, ballistics \cite{Hall88}, topological data analysis \cite{Bobrowski17,Tillmann22}, stochastic optimization \cite{ZZ}, immunology \cite{Moran}, wireless communications \cite{Lan}, ground states of interacting particles \cite{Torquato10}, and other subjects. 

In this work, we explore {\em random} sequential coverings (RSCs). Previous studies of coverings, particularly coverings of the ring by randomly placed arcs \cite{Stevens39,Domb47,Dvoretzky56,Flatto,Mandelbrot72,Jonasson,Jonasson08,hitting}, focused on the probability of the complete covering, i.e., the dependence on the average coverage from the number of intervals. In contrast, we consider the RCS processes dynamically: The objects are deposited sequentially until the complete coverage is reached~\footnote{Very different dynamics was studied in Refs.~\cite{Jonasson,Jonasson08}.}. The chief difference with previous work is that a deposition attempt is accepted only if it leads to an increase in coverage. 

One can study the RSCs by objects of arbitrary shape, size and orientation. The RSCs of the plane by discs or aligned squares exemplify coverings by objects of fixed shape, size and orientation; the RSCs of the plane by ellipses or squares of arbitrary orientation are examples of coverings by identical objects of arbitrary orientation. We limit ourselves with identical objects of fixed orientation. Coverings by objects of varying size are particularly natural in one dimension, see \cite{Shepp72a,Shepp72b,Siegel78,Siegel79,Holst,Shepp87} studying coverings of the ring by arcs with random length. 

We mostly focus on one dimension where the RCS processes are analytically tractable. The one-dimensional RCS processes shed light on the properties of the RCS processes in higher dimensions. We begin with discrete coverings of finite intervals. An interval of size $L$ is a set of $L$ consecutive integers, deposited identical objects could be dimers $\{i,i+1\}$, or trimers $\{i,i+1,i+2\}$, generally $\ell-$mers. In Sect.~\ref{sec:Dimers}, we consider the RSC of a finite interval by dimers. We compute the average total number of dimers in a congested configuration using the same methods as in the analysis of the random sequential packing of the dimers. In Sec.~\ref{sec:Dimers-Fluct}, we probe fluctuations in the total number of dimers in congested configurations adopting again an approach used in the packing problem \cite{Krapivsky20}. In particular, we compute the full counting statistics and extreme probabilities of maximally and minimally congested configurations. In Sec.~\ref{sec:L-mers}, we study the RSC of a finite interval by segments with $\ell$ sites, $\ell-$mers, for arbitrary integer $\ell \geq 2$. 

In Sec.~\ref{sec:inf}, we study the jamming properties of the covering of the infinite lattice $\mathbb{Z}$ by $\ell-$mers. In Sec.~\ref{sec:inf-time}, we analyze the dynamics of the covering of the one-dimensional lattice by $\ell-$mers. We show that the fraction $\pi_0(t)$ of uncovered sites exhibits a pure exponential decay: $\pi_0= e^{-\ell t}$. We also determine the average coverage accounting for multiplicity, $\sum_{0\leq k\leq \ell} k\pi_k(t)$, where $\pi_k(t)$ is the fraction of sites covered $k$ times. 

For the RSC by dimers we compute (Sec.~\ref{subsec:dimers}) all fractions $\pi_0, \pi_1, \pi_2$. For the RSC of the one-dimensional lattice by $\ell-$mers, the non-trivial fractions $\pi_1, \ldots, \pi_\ell$ are unknown when $\ell\geq 3$. In Sec.~\ref{subsec:trimers} and Appendix~\ref{app:small}, we present conjecturally exact results for the fractions $\pi_1, \pi_2, \pi_3$ for the RSC by trimers ($\ell=3$). The challenge of computing  $\pi_1, \ldots, \pi_\ell$ suggests
modifying the deposition rule, namely, accepting only those deposition events where the overlap with previously deposited $\ell-$mers is at most $\lfloor \frac{\ell}{2}\rfloor$. This model B is more tractable since for any $\ell$, only $\pi_0, \pi_1, \pi_2$ are non-trivial (Sec.~\ref{sec:mod-B}).

In Sec.~\ref{sec:line}, we turn to the continuous RSC of the line by sticks of unit length. The fractions $\pi_k(t)$ of line covered $k$ times are non-trivial for all $k\geq 0$. We compute the fraction $\pi_0(t)$ of the uncovered line and the average coverage $\sum_{k\geq 1} k\pi_k(t)$ accounting for multiplicity; the fractions $\pi_k$ with $k\geq 1$ are still unknown. We also examine the continuous RSC with model B dynamics in which a deposition attempt is accepted only when the center of the incoming stick lands in the uncovered region. In this model, each point is covered by at most two sticks. We compute $\pi_0, \pi_1, \pi_2$.  In Sec.~\ref{sec:space}, we analyze the continuous RSCs of $\mathbb{R}^d$ by identical balls. We establish the asymptotic decay of the uncovered volume fraction $\pi_0(t)$.  In Sec.~\ref{sec:CR}, we briefly discuss related models and challenges for future work.

\section{Coverings by Dimers: Average Properties}
\label{sec:Dimers}

Packing the one-dimensional lattice by dimers placed at random without overlaps has a long history going back to Flory \cite{Flory39}, see \cite{Evans93,Talbot00,KRB} for review of this and other random sequential adsorption processes. 

In this section, we analyze the RSC of an interval by dimers. We call a set of $L$ consecutive integers an interval of size $L$. We shortly write $[1,L]$ for the interval $\{1,\ldots,L\}$.  An interval of size two is a dimer. In a successful deposition event at least one previously uncovered site of the interval is covered; unsuccessful events are discarded.For instance, successful sequential (top to bottom) deposition events 
\begin{equation}
\label{illustr}
\begin{split}
& \bullet\bullet\,\circ\circ\,\circ\\
& \bullet\bullet\,\circ\circ\bullet\bullet\\
& \bullet\bullet\,\circ\bullet\,\blacktriangle\,\bullet\\
& \bullet\blacktriangle\bullet\bullet\,\blacktriangle\,\bullet
\end{split}
\end{equation}
illustrate the RSC of an empty interval $[1,5]$ by dimers: The first dimer is $[1,2]$, the second  $[5,6]$, the third $[4,5]$, and the fourth $[2,3]$. An empty site is denoted by $\circ$, a site covered by one dimer is denoted by $\bullet$, and a site covered by two dimers is denoted by $\blacktriangle$. 

Generally, we start with an empty interval $[1,L]$. The first dimer $[k-1,k]$ with $k$ chosen randomly among the sites $k=1,\ldots,L+1$. The second successfully deposited dimer $[j-1,j]$ with $j$ randomly chosen among $1,\ldots,L+1$ with $j\ne k$ due to the requirement that the coverage must increase. The third successfully deposited dimer $[i-1,i]$ with $i$ randomly chosen among $1,\ldots,L+1$ and satisfying $i\ne j$ and $i \ne k$; in the special case that $j=k-2$, the requirement that the coverage increases also implies that $i\ne k-1$, and similarly if $j=k+2$, it should be $i\ne k+1$. This procedure continues until a congested configuration is reached, that is, all sites $1,\ldots,L$ become covered. 

In the final congested state all sites $[1,\ldots,L]$ are covered. Sites $0$ and $L+1$  can also be covered; if this happens we have one or two overhangs. The congested configuration in \eqref{illustr} has a single overhang on the right. 

\subsection{The average number of dimers}

Denote by $N$ the total number of dimers in a congested configuration; equivalently, $N$ is the total number of (successful) deposition events. When $L=1$, one deposition event suffices: $N=1$. For $L\geq 2$, the total number of deposition events varies from realization to realization. For instance, $N=1$ (probability $1/3$) or $N=2$ (probability $2/3$) if $L=2$. When $L=3$, one gets $N=2$ (probability $2/3$) or  $N=3$ (probability $1/3$). The bounds on $N$ are 
\begin{equation}
\label{bounds}
\Big\lfloor\frac{L+1}{2}\Big\rfloor \leq N \leq L
\end{equation}
where $\lfloor x \rfloor$ denotes the largest integer not exceeding $x$. The upper bound in \eqref{bounds} follows from the fact that in each successful deposition event at least one previously uncovered site gets covered. The maximal number of added covered sites in a successful deposition event is two leading to the lower bound in \eqref{bounds}. 

Let us compute the average number $D_L=\langle N\rangle$ of deposition events before the interval gets congested. We already know $D_1=1, ~D_2 = \frac{5}{3}, ~D_3 = \frac{7}{3}$. Suppose the first dimer is $[k-1,k]$ with $k$ uniformly chosen in the range $1\leq k\leq L+1$. Subsequent coverings of the intervals $[1,k-2]$ and $[k+1,L]$ on the left and right of the first dimer proceed independently. This key feature makes the one-dimensional process tractable. 

In the case of $D_L$, we arrive at the recurrence
\begin{equation}
\label{DL:rec}
D_L = \frac{1}{L+1}\sum_{k=1}^{L+1} \big(D_{k-2}+1+D_{L-k}\big)
\end{equation}
The $(L+1)^{-1}$ factor is the probability that the first dimer is $[k-1,k]$. In $D_{k-2}+1+D_{L-k}$, the first term $D_{k-2}$ accounts dimers covering $[1,k-2]$, to which we add one for the first dimer and $D_{L-k}$ accounting for the average number of dimers covering $[k+1,L]$. Summing over all possible $k=1,\ldots,L+1$ gives the recurrence \eqref{DL:rec}. 

One verifies that \eqref{DL:rec} is applicable for all $L\geq 1$ if we set $D_{-1}=D_0=0$. The generating function 
\begin{equation}
\label{Dx:def}
D(x) = \sum_{L\geq 1} D_{L}\, x^{L+1}
\end{equation}
allows us to convert the recurrence \eqref{DL:rec} into a differential equation. Indeed, multiplying \eqref{DL:rec} by $(L+1)x^L$ and summing over all $L\geq 1$ one reduces the left-hand side to $\frac{dD}{dx}$. The right-hand side also simplifies and one obtains 
\begin{equation}
\frac{dD}{dx} = \frac{2}{1-x}\,D+(1-x)^{-2}-1
\end{equation}
Solving this linear inhomogeneous differential equation subject to the initial condition $D(0)=0$ yields
\begin{equation}
\label{Dx:sol}
D(x) = \frac{2}{3}\,\frac{x^3}{(1-x)^2} + \frac{x^2}{1-x}
\end{equation}
Comparing \eqref{Dx:sol} with the definition \eqref{Dx:def} we extract 
\begin{equation}
\label{DL:sol}
D_L = \frac{2}{3}\,L + \frac{1}{3}
\end{equation}
for all $L\geq 1$. 

\subsection{Behavior near the boundary}

Let us compute the fraction of the coverings without overhang on the left. Equivalently, we want to  compute  the probability $p_L$ that $[1,2]$ is the left-most dimer. The probabilities $p_j$ satisfy the recurrence
\begin{equation}
\label{pL:rec}
p_L =  \frac{1}{L+1} + \frac{1}{L+1}\sum_{k=2}^{L} p_{k-1}
\end{equation}
Indeed, $(1,2)$ may be the first deposited dimer. This event occurs with probability $(L+1)^{-1}$ and explains the first term on the right-hand side of \eqref{pL:rec}. If $[k,k+1]$ with $2\leq k\leq L$  is the first deposited dimer, $[1,2]$ will be the left-most dimer with probability $p_{k-1}$. This explains the sum on the right-hand side of \eqref{pL:rec}. The recurrence \eqref{pL:rec} admits a simple solution: 
\begin{equation}
\label{pL:sol}
p_L = \frac{1}{2}
\end{equation}
The fraction of the coverings without overhang on the right is also $p_L$. 

Let $q_L$ be the fraction of the coverings without overhangs. For small $L$, one can compute $q_L$ by hand to yield 
\begin{equation}
q_1=0,\quad q_2=\frac{1}{3}, \quad q_3 = \frac{1}{4}
\end{equation}
etc. Generally for $L\geq 3$, the probability $q_L$ satisfies the recurrence
\begin{equation}
\label{qL:rec}
q_L =  \frac{2}{L+1}\,p_{L-2} + \frac{1}{L+1}\sum_{k=2}^{L-2} p_{k-1}p_{L-k-1}
\end{equation}
which is established similarly to the recurrence \eqref{pL:rec}. By inserting \eqref{pL:sol} into \eqref{qL:rec} we find that $q_L$ stabilizes for $L\geq 3$:
\begin{equation}
\label{qL:sol}
q_L = \frac{1}{4}
\end{equation}

\subsection{The total number of congested configurations}

 In the RSC process, different congested configurations occur with different probabilities, but here we ignore this feature and just determine the total number of congested dimer configurations $C_L$. Congested configurations can be divided into two complementary sets: (i) congested configurations with the leftmost $[0,1]$ dimer, (ii) congested configurations with the leftmost $[1,2]$ dimer. In the former case, we need to clog the interval $[2,L]$;  in the latter case, the interval $[3,L]$ must get clogged. This leads to the recurrence
\begin{equation}
\label{CL-rec}
C_L = C_{L-1}+C_{L-2}
\end{equation}
defining Fibonacci numbers. The initial conditions are $C_1=2$ and $C_2=4$, and the solution of \eqref{CL-rec} satisfying these initial conditions is
\begin{equation}
\label{CL-Fib}
C_L = 2 F_{L+1}
\end{equation}
where $F_n$ are standard Fibonacci numbers.

\subsection{Coverage of rings by dimers}

We consider the RSC of a lattice interval $[1,L]$ if not mentioned otherwise, but in this subsection we briefly discuss the RSC of a ring with $L$ sites. This process can be studied using the same approaches, and some answers can be extracted from the previous results for the RSC of intervals. Here we present a few results in the case of the dimer RSC of rings. 

We set $L\geq 2$ to allow deposition of dimers. The number of dimers in a congested configuration, equivalently the total number of successful deposition events required to reach a congested covering of a ring, is deterministic when $L=2$ and $L=3$, viz. $N=1$ when $L=2$ and $N=2$ when $L=3$.  For $L\geq 4$, the total number of successful deposition events varies from realization to realization. 

The average number $N_L$ of dimers in a congested coverage of the ring with $L$ sites can be expressed through the average number $D_{L-2}$ of dimers in a congested coverage of the interval with $L-2$ sites that are uncovered after the first deposition  event. One gets $N_L=1+D_{L-2}$ which in conjunction  with \eqref{DL:sol} yields
\begin{equation}
\label{DL:ring}
N_L = \frac{2}{3}\,L
\end{equation}
for $L\geq 3$. The bounds on $N$ in the case of a ring read
\begin{equation}
\label{bounds-ring}
\Big\lfloor\frac{L+1}{2}\Big\rfloor \leq N \leq L-1
\end{equation}

The version of the RSC process in which all attempts are accepted is popular in mathematics literature. Suppose we randomly place $\mathcal{N}$ dimers on the ring.  One basic quantity is the probability $P(\mathcal{N},L)$ that after $\mathcal{N}$ deposition events the ring is fully covered. This probability is non-trivial, $0<P(\mathcal{N},L)<1$, for all
\begin{equation}
\label{bounds-all}
\Big\lfloor\frac{L+1}{2}\Big\rfloor \leq \mathcal{N} < \infty
\end{equation}

Another basic quantity is the average fraction $\pi(\mathcal{N},L)$ of covered sites. The quantities $P(\mathcal{N},L)$ and $\pi(\mathcal{N},L)$ have been studied \cite{Stevens39,Domb47,Flatto,Siegel78,Siegel79,Holst,Shepp87} in the continuous case when intervals of unit length are randomly placed into a ring of an arbitrary length $L$. In particular, in the continuous case the probability $P(\mathcal{N},L)$ has been explicitly evaluated already in the classical paper by Stevens \cite{Stevens39}. In our situation when the ring is covered by dimers, the probability that $\mathcal{N}$ deposition events do not cover some arbitrary site is $\left(\frac{L-2}{L}\right)^\mathcal{N}$, and  therefore
\begin{equation}
\label{ring-covered}
\pi(\mathcal{N},L) =  1 - \left(\frac{L-2}{L}\right)^\mathcal{N}
\end{equation}

\section{Coverings by Dimers: Fluctuations}
\label{sec:Dimers-Fluct}

One-dimensional RSC models are sufficiently simple, and in addition to the average characteristics one can analytically probe fluctuations. 

\subsection{Full counting statistics}

The total number $N$ of deposited dimers fluctuates from one congested configuration to another. The cumulant generating function encodes all the cumulants of $N$. By definition, the cumulant generating function is 
\begin{equation}
\label{F:def}
F(\lambda, L) \equiv \langle e^{\lambda N}\rangle = \sum_N e^{\lambda N} P(N,L)
\end{equation}
where $P(N,L)$ is the probability to have $N$ dimers in the congested configuration. The standard relation
\begin{equation}
\ln \langle e^{\lambda N}\rangle = \sum_{n\geq 1} \frac{\lambda^n}{n!}\,  \langle N^n\rangle_c
\end{equation}
then gives all the cumulants: the average $\langle N\rangle_c = \langle N\rangle$, the variance $\langle N^2\rangle_c = \langle N^2\rangle - \langle N\rangle^2$, etc. 

The function $F(\lambda, L) \equiv \langle e^{\lambda N}\rangle$ grows as $ e^{LU(\lambda)}$ when $L\gg 1$. More precisely 
\begin{equation}
\label{U:def}
U(\lambda) = \lim_{L\to\infty} \frac{\ln F(\lambda, L)}{L}
\end{equation}
Expanding $U(\lambda)$ into the Taylor series near $\lambda=0$ we can read off the cumulants: 
\begin{equation}
\label{U-exp}
U(\lambda)=\sum_{n\geq 1} \frac{\lambda^n}{n!}\, U_n, \qquad \langle N^n\rangle_c = L U_n
\end{equation}

To determine $F(\lambda, L)$ we proceed in the same way as in the derivation of the recurrence \eqref{DL:rec} for the  average. The recurrence 
\begin{equation}
\label{FL:rec}
F(\lambda, L) = \frac{e^\lambda}{L+1}\sum_{k=1}^{L+1} F(\lambda, k-2) F(\lambda, L-k)
\end{equation}
is established by noticing that after the first deposition event, viz. putting a dimer at $(k-1,k)$, the intervals on the left and on the right are filled independently. 

\begin{figure}[ht]
\begin{center}
\includegraphics[width=0.44\textwidth]{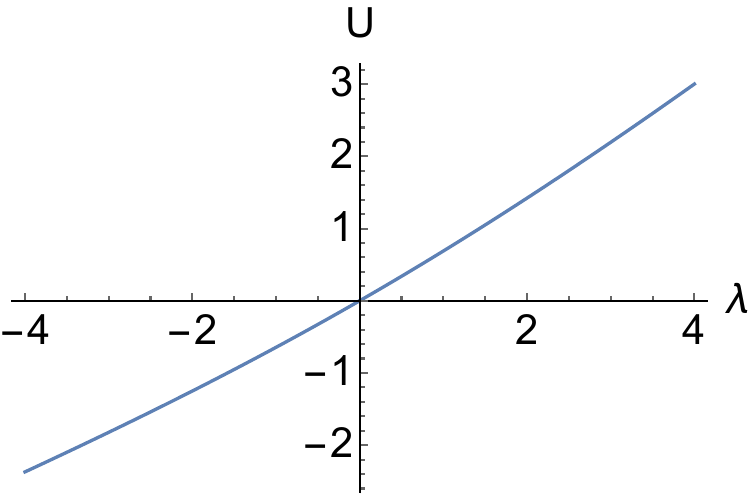}
\caption{The plot of the cumulant generating function $U(\lambda)$. The expansion of $U(\lambda)$ at $\lambda=0$ yields the cumulants.}
\label{fig:UL_cov}
  \end{center}
\end{figure}

The recurrence \eqref{FL:rec} applies for all $L\geq 1$ if we set
\begin{equation}
\label{IC:FF}
F(\lambda, -1) = F(\lambda, 0) = 1
\end{equation}
We now introduce the generating function
\begin{equation}
\label{Phi:def}
\Phi(\lambda, x) = \sum_{L\geq 0}  F(\lambda, L)\, x^{L+1}
\end{equation}
To recast \eqref{FL:rec} into an equation for the generating function we multiply \eqref{FL:rec} by $(L+1)x^{L}$ and sum over all $L\geq 1$. After a bit of algebra we arrive at a differential equation
\begin{equation}
\label{Phi:eq}
\frac{\partial \Phi}{\partial x} = 1 + e^\lambda \Phi(2+\Phi)
\end{equation}
Solving  \eqref{Phi:eq} we get
\begin{equation}
\label{Phi:sol}
\Phi(\lambda, x) = \frac{\tan(\Lambda x)}{\Lambda-e^{\lambda}\tan(\Lambda x)}\,, \quad \Lambda\equiv e^{\lambda/2}\sqrt{1-e^\lambda}
\end{equation}
The generating function $\Phi(\lambda, x)$ has a simple pole at 
\begin{equation}
\label{y:def}
x=y(\lambda) = \frac{\arctan\big(e^{-\lambda/2}\sqrt{1-e^\lambda}\big)}{e^{\lambda/2}\sqrt{1-e^\lambda}}
\end{equation}
The cumulant generating function 
\begin{equation}
\label{U:implicit}
U(\lambda) = -\ln y(\lambda)
\end{equation}
is plotted in Fig.~\ref{fig:UL_cov}. 

Expanding $U(\lambda)$ in powers of $\lambda$, cf. Eq.~\eqref{U-exp}, gives the already known value $\frac{\langle N\rangle}{L} = \frac{2}{3}$ in the leading order.  The variance and the next two cumulants read
\begin{equation*}
\frac{\langle N^2\rangle_c}{L}  = \frac{2}{45}\,, \quad 
\frac{\langle N^3\rangle_c}{L}  = \frac{2}{945}\,, \quad
\frac{\langle N^4\rangle_c}{L}  =  -\frac{22}{4725}
\end{equation*}

The ratios $\langle N^n\rangle_c/\langle N\rangle$ of cumulants to the average are known as Fano factors \cite{Fano}.
Here are a few Fano factors 
\begin{equation}
\label{Fano}
\begin{split}
\frac{\langle N^2\rangle_c}{\langle N\rangle} & = \frac{1}{15}\\
\frac{\langle N^3\rangle_c}{\langle N\rangle} & = \frac{2}{315}\\
\frac{\langle N^4\rangle_c}{\langle N\rangle} & =  -\frac{11}{1575}\\
\frac{\langle N^5\rangle_c}{\langle N\rangle} & =  -\frac{1}{1485}\\
\frac{\langle N^6\rangle_c}{\langle N\rangle} & =  \frac{47221}{14189175}\\
\frac{\langle N^7\rangle_c}{\langle N\rangle} & =  \frac{811}{2027025}\\
\frac{\langle N^8\rangle_c}{\langle N\rangle} & =  -\frac{1790851}{516891375}
\end{split}
\end{equation}
If $N$ were a Poisson distributed random variable, all cumulants would be equal,  $\langle N^n\rangle_c = \langle N\rangle$, and all Fano factors equal to unity. 

The so-called  Mandel $Q$ parameter \cite{Mandel79} defined via
\begin{equation}
\label{Q:def}
Q=\frac{\langle N^2\rangle_c}{\langle N\rangle}-1
\end{equation}
is a basic measure characterizing the deviation from Poissonian statistics. The values  $-1\leq Q<\infty$ are permissible. For the Poisson statistics $Q=0$ and the range $-1\leq Q<0$ is sub-Poissonian. We have
\begin{equation}
\label{Q:RSA}
Q = -\frac{14}{15}
\end{equation}
indicating that the statistics of dimer covering is strongly sub-Poissonian. 

\subsection{Extremal congested configurations}
\label{subsec:Extrem}

Consider extremal congested configurations. The minimal number of dimers covering the interval is $\lfloor\frac{L+1}{2}\rfloor$, see \eqref{bounds}. For the interval with even number of sites, say $L=2n$, we have $N_\text{min}=n$. The probability $m_n$ of reaching this configuration can be determined from recurrence 
\begin{equation}
\label{m-n}
m_n = \frac{1}{2n+1}\sum_{k=0}^{n-1}m_k m_{n-k-1}
\end{equation}
The boundary condition is $m_0 =  1$. 

Introducing the generating function 
\begin{equation}
\label{m:GF}
m(x)=\sum_{k\geq 0}m_k x^{k}
\end{equation}
we convert the recurrence \eqref{m-n} into a Riccati equation 
\begin{equation}
\label{m:Ric}
2\,\frac{dm}{dx}=\frac{1-m}{x}+m^2
\end{equation}
Solving \eqref{m:Ric} subject to $m(0)=1$ we obtain
\begin{equation}
\label{m:GF-sol}
m(x) = \frac{\tan \sqrt{x}}{\sqrt{x}}
\end{equation}
Analyzing the divergence of the generating function when $x\to (\pi/2)^2$, we extract the large $n$ asymptotic
\begin{equation}
\label{mu-asymp}
m_n \simeq 2\left(\frac{2}{\pi}\right)^{2n+2}
\end{equation}

The maximal number of dimers covering the interval with $L$ sites is $N_\text{max}=L$. For such covering to arise, the first deposition event must be on the edge of the interval covering just one site. The following events must be on the edges of the remaining interval, and this kind of 
avalanche must proceed. The probability is 
\begin{equation}
\label{M-asymp}
M_L = \prod_{j=2}^{L+1} \frac{2}{j} = \frac{2^L}{(L+1)!}
\end{equation}

The asymptotic behaviors of the extreme probabilities of maximally and minimally congested configurations can be also extracted from  the asymptotic behavior of the cumulant generating function in the $\lambda\to \pm \infty$ limits. We merely outline this method as the direct derivation is simpler and gives more precise results \eqref{mu-asymp} and \eqref{M-asymp}. An asymptotic analysis of \eqref{y:def} yields
\begin{equation}
\label{U:asympt}
U(\lambda) = 
\begin{cases}
\frac{1}{2}\lambda - \ln(\pi/2)        & \lambda\to -\infty\\
\lambda - \ln(\ln\lambda)               & \lambda\to \infty
\end{cases}
\end{equation}
where we have dropped terms vanishing in the $\lambda\to \pm \infty$ limits.  Using \eqref{F:def} and \eqref{U:def} together with the asymptotic of $U(\lambda)$ in the $\lambda\to - \infty$ limit, one recovers \eqref{mu-asymp}. Similarly in the  $\lambda\to \infty$ limit, one can recover \eqref{M-asymp}.

\section{Coverings by $\ell-$mers}
\label{sec:L-mers}

Let us generalize the deposition of dimers to the deposition of sub-intervals with $\ell$ sites, $\ell-$mers. Instead of Eq.~\eqref{DL:rec} we obtain
\begin{equation}
\label{DLL:rec}
D_L = \frac{1}{L+\ell-1}\sum_{k=1}^{L+\ell -1} \big(D_{k-\ell}+1+D_{L-k}\big)
\end{equation}
Using the generating function 
\begin{equation}
\label{DxL:def}
D(x) = \sum_{L\geq 1} D_{L}\, x^{L+\ell - 1}
\end{equation}
we recast the recurrence \eqref{DLL:rec} into
\begin{equation}
\frac{dD(x)}{dx} = \frac{2}{1-x}\,D(x)+(1-x)^{-2}-\sum_{j=1}^{\ell - 1}jx^{j-1}
\end{equation}
from which
\begin{equation}
D(x) = \frac{2}{\ell + 1}\,\frac{x^{\ell + 1}}{(1-x)^2} + \frac{x^{\ell}}{1-x}
\end{equation}
Thus
\begin{equation}
\label{DLL:sol}
D_L = \frac{2}{\ell + 1}\,L + \frac{\ell-1}{\ell + 1}
\end{equation}

The bounds \eqref{bounds} generalize to
\begin{equation}
\label{bounds:L}
\Big\lfloor\frac{L+\ell-1}{\ell}\Big\rfloor \leq N \leq L
\end{equation}
In particular, the minimal number of $\ell-$mers covering the interval with $L=n\ell$ sites is $N_\text{min}=n$. The probability $m_n$ of reaching this configuration can be determined from recurrence 
\begin{equation}
\label{m-nL}
m_n = \frac{1}{\ell n+\ell-1}\sum_{k=0}^{n-1}m_k m_{n-k-1}
\end{equation}
with boundary condition $m_0 =  1$. The generating function \eqref{m:GF} satisfies a Riccati equation 
\begin{equation}
\label{mL:Ric}
\ell\,\frac{dm}{dx}=(\ell-1)\frac{1-m}{x}+m^2
\end{equation}
This Riccati equation admits an explicit solution through the Bessel function:
\begin{equation}
\label{m-L:sol}
m = -\ell\,\frac{d}{dx} \,\ln\!\left[J_{-1/\ell}\!\left(\frac{2\sqrt{(\ell-1)x}}{\ell}\right)\right] - \frac{1}{2x}
\end{equation}
The generating function \eqref{m-L:sol} diverges when $x\to y$, where $y=y(\ell)$ is found from 
\begin{equation}
J_{-1/\ell}\!\left(\frac{2\sqrt{(\ell-1)y(\ell)}}{\ell}\right) = 0
\end{equation}
One can verify that the generating function has a simple pole at $x=y(\ell)$ and hence $m_n  \sim [y(\ell)]^{-n}$. Re-writing in terms of the length of the chain, we arrive at 
\begin{equation}
\label{m-y}
 P(N_\text{min},L)\sim [u(\ell)]^{-L}, \qquad u(\ell)=[y(\ell)]^{1/\ell}
\end{equation}

\begin{figure}[ht]
\begin{center}
\includegraphics[width=0.456\textwidth]{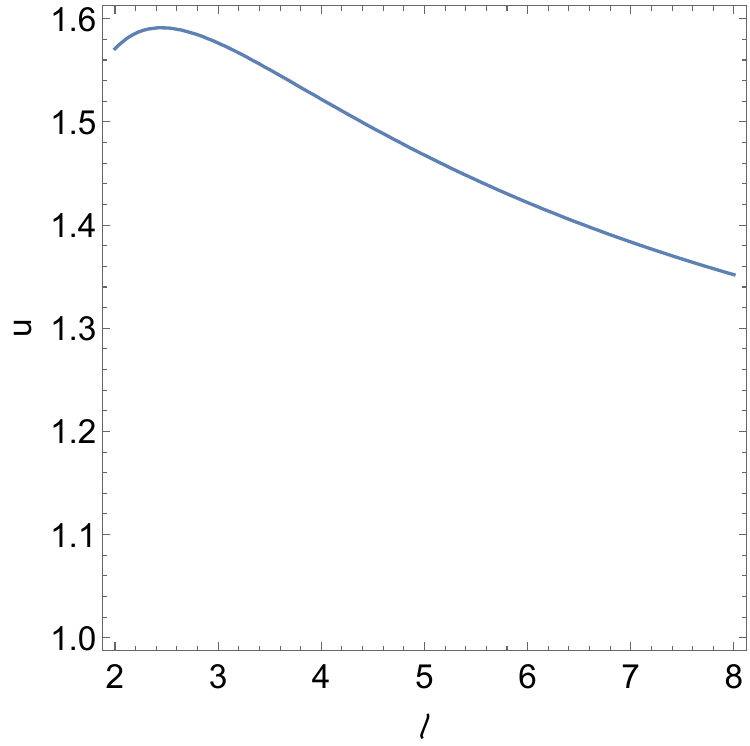}
\caption{The plot of $u(\ell)$ appearing in Eq.~\eqref{m-y}  giving the probability to reach the congested configuration with minimal number of $\ell-$mers. Only integer values matters; in the present plot, $\ell=2, \ldots, 8$.}
\label{fig:u28}
  \end{center}
\end{figure}

The plot of $u(\ell)$ is drawn [Fig.~\ref{fig:u28}] for all $\ell\geq 2$, albeit only positive integer values, $\ell\in \mathbb{Z}_+$, matter. A few first numerical values are
\begin{equation*}
\begin{split}
u(2) & = 1.57079632679\ldots\\
u(3) & = 1.57656918868\ldots\\
u(4) & = 1.52206284731\ldots\\
u(5) & = 1.46785551917\ldots\\
u(6) & = 1.42179975053\ldots\\
u(7) & = 1.38361356892\ldots
\end{split}
\end{equation*}
We also note a large $\ell$ asymptotic:
\begin{equation}
\label{uL:large}
u(\ell)\simeq \left[\frac{a_1}{2}\,\frac{\ell}{\sqrt{\ell-1}}\right]^{2/\ell}
\end{equation}
where $a_1=2.40482555769577\ldots$ is the first positive zero of the Bessel function $J_0(a_1)=0$. Using \eqref{uL:large} one finds 
\begin{equation}
u(\ell) - 1 \simeq \ell^{-1}\,\ln\!\left[\frac{a_1^2}{4}\,\ell\right]
\end{equation}
when $\ell\gg 1$.

The maximal number of $\ell-$mers covering the interval with $L$ sites is $N_\text{max}=L$ and the probability to reach such jammed configurations is
\begin{equation}
M_L = \prod_{j=\ell}^{L+\ell-1} \frac{2}{j} = \frac{2^L(\ell-1)!}{(L+\ell-1)!}
\end{equation}

\section{Infinite Lattice}
\label{sec:inf}

When $L=\infty$, the average number of $\ell-$mers covering any point is
\begin{equation}
\lim_{L\to \infty}\frac{\ell D_L}{L} = \frac{2\ell}{\ell+1}
\end{equation}
Denote by $\pi_k^{(\ell)}$ the fraction of sites  with coverage $k$. We have $\pi_k^{(\ell)}=0$ when $k=0$ and $k>\ell$. The sum rule
\begin{subequations}
\begin{equation}
\label{norm}
\sum_{k=1}^\ell \pi_k^{(\ell)} = 1
\end{equation}
reflects the normalization. Another sum rule
\begin{equation}
\label{cov}
\sum_{k=1}^\ell k\,\pi_k^{(\ell)} =  \frac{2\ell}{\ell+1}
\end{equation}
\end{subequations}
expresses the average coverage. When $\ell=2$, these sum rules fix the coverage distribution
\begin{equation}
\label{coverage-dimers}
\pi_1^{(2)} = \frac{2}{3}\,, \qquad \pi_2^{(2)} = \frac{1}{3}
\end{equation}

The fractions $\pi_k^{(\ell)}$ are unknown when $\ell\geq 3$.  In Sec.~\ref{sec:inf-time} we {\em guess} the temporal behavior of the fractions for $\ell=3$ from which
\begin{equation}
\label{coverage-trimers}
\pi_1^{(3)} = \frac{2}{3}\,, \qquad \pi_2^{(3)} = \pi_3^{(3)} = \frac{1}{6}
\end{equation}
for the RSC of the one-dimensional lattice $\mathbb{Z}$ by trimers. 

Thus we know $\pi_k^{(\ell)}$ only for dimer coverage when single and double coverages are the only possibilities. This feature hints that one-dimensional RSC processes could be more tractable if {\em only} single and double coverage were permitted. We call the original RSC process model A and define model B by postulating that a deposition event is accepted only when overlap with previously deposited $\ell-$mers is at most $\lfloor \frac{\ell}{2}\rfloor$. When $\ell$ is odd, an equivalent definition is that the center of a new $\ell-$mer lands into an uncovered site. In the realm of model B, any site  can be only single and double-covered. 

In the case of the dimer coverage, model B is identical to model A, so the coverage distribution is given by \eqref{coverage-dimers}. In the case of the trimer coverage, the temporal behavior of the model B (see Sec.~\ref{subsec:trimers-B}) leads to the final fractions
\begin{equation}
\label{cov-3-B}
\pi_1^{(3)} = \frac{1+3e^{-2}}{2}\,, \qquad  \pi_2^{(3)} = \frac{1-3e^{-2}}{2}
\end{equation}

The same approach (see Sec.~\ref{subsec:trimers-B}) in principle allows one to determine the fractions $\pi_1^{(\ell)}$ and  $\pi_2^{(\ell)}$ for the model B dynamics when $\ell$ is arbitrary. The calculations are straightforward but become unwieldy as $\ell$ increases. We thus present the fractions only for $\ell=4$ and $\ell=5$. For the RSC by $4-$mers
\begin{equation}
\label{cov-4-B}
\pi_1^{(4)} = 1-3e^{-2}\,, \qquad  \pi_2^{(4)} = 3e^{-2}
\end{equation}
For the RSC by $5-$mers
\begin{equation}
\label{cov-5-B}
\pi_2^{(5)} = \frac{5\sqrt{\pi}}{2\,e^4}\,\{\text{Erfi}[2]-\text{Erfi}[1]\} - 1 = 0.3727549\ldots
\end{equation}
where
\begin{equation}
\label{error}
\text{Erfi}[z] = \frac{2}{\sqrt{\pi}}\int_0^z dy\,e^{y^2}
\end{equation}
is an imaginary error function.

\section{Temporal Behavior}
\label{sec:inf-time}

In the preceding sections, we have studied congested configurations. The covering algorithm is dynamical, so here we analyze the evolution. A finite interval gets congested in a finite (albeit random) time. If the substrate is infinite, $\mathbb{Z}$ or $\mathbb{R}$ in one dimension, the RSC processes continue indefinitely. In this section we study the dynamics of RSC of the one-dimensional lattice. 

\subsection{Coverings by dimers}
\label{subsec:dimers}

Let $E_m(t)$ be the density of empty strings of length $m$:
\begin{equation}
\label{Em}
E_m = \text{Prob}[\underbrace{\circ \cdots\circ}_m]
\end{equation}
(Recall that $\circ$ denotes an empty site.) In the case of the RSC by dimers, $E_m(t)$ satisfies the master equation
\begin{equation}
\label{Emt}
\dot E_m =-(m+1)E_m, \qquad m\geq 1
\end{equation}
where the over-dot denotes time derivative and we have set the deposition rate to unity. Solving \eqref{Emt} subject to the initial condition $E_m(0)=1$ gives \begin{equation}
\label{Em-sol}
E_m(t)=e^{-(m+1)t}
\end{equation}
In particular, the density $\pi_0(t) $ of empty sites is 
\begin{equation}
\label{E-sol}
\pi_0(t) \equiv E_1(t) = e^{-2t}
\end{equation}

Let $V_m$ be the density of voids of length $m$ 
\begin{equation}
\label{Vm}
V_m = \text{Prob}[\bullet\underbrace{\circ \cdots\circ}_m \bullet]
\end{equation}
where $\bullet$ denotes a single-covered site. Expressing the density of voids through the density of empty strings (see, e.g., \cite{KRB}) we obtain
\begin{equation}
\label{V-E}
V_m = E_m - 2E_{m+1} + E_{m+2}
\end{equation}
Combining \eqref{V-E} with \eqref{Em-sol} we get 
\begin{equation}
\label{Vm-sol}
V_m(t)=e^{-(m+1)t}(1-e^{-t})^2
\end{equation}

The density $\pi_1(t)$ of single-covered sites grows according to the evolution equation
\begin{equation}
\label{S-eq}
\frac{d \pi_1}{dt}  = 2 \sum_{m\geq 1} (m-1) V_m 
\end{equation}
Indeed, there are $m+1$ deposition events destroying the void \eqref{Vm}. The boundary sites in \eqref{Vm} are single-covered, so two deposition events do not change the number of single-covered sites, while each of the remaining  $m-1$ deposition events increase the number of single-covered sites by two. Plugging \eqref{Vm-sol} into \eqref{S-eq} and integrating we obtain $\dot \pi_1 = 2 e^{-3t}$ from which
\begin{equation}
\label{S-sol}
\pi_1(t) = \tfrac{2}{3} - \tfrac{2}{3} e^{-3t}
\end{equation}
The asymptotic value $\pi_1(t=\infty)=\frac{2}{3}$ agrees with \eqref{coverage-dimers}.

The density $\pi_2(t)$ of double-covered sites varies according to the evolution equation
\begin{equation}
\label{D-eq}
\frac{d  \pi_2}{dt} = 2 \sum_{m\geq 1} V_m
\end{equation}
Indeed, the boundary sites in \eqref{Vm} are single-covered, and each of them can become double-covered after deposition events destroying the void \eqref{Vm} and covering this  site. Plugging \eqref{Vm-sol} into \eqref{D-eq} and integrating we obtain
\begin{equation}
\label{D-sol}
\pi_2(t) = \tfrac{1}{3}  - e^{-2t} + \tfrac{2}{3} e^{-3t}
\end{equation}
The asymptotic value $\pi_2(t=\infty)=\frac{1}{3}$ agrees with \eqref{coverage-dimers}. 

\begin{figure}[ht]
\begin{center}
\includegraphics[width=0.456\textwidth]{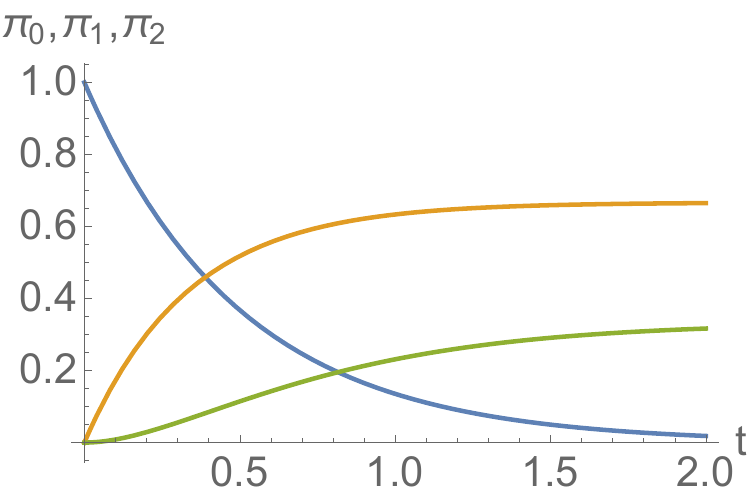}
\caption{The RSC by dimers. Shown is the evolution  of the densities $\pi_0, \pi_1, \pi_2$ of uncovered, single-covered, and double-covered sites. The densities $\pi_1$ and $\pi_2$ are increasing functions of time satisfying $\pi_1(t)>\pi_2(t)$ for all $t>0$. The density $\pi_0(t)$ decays exponentially, Eq.~\eqref{E-sol}. }
\label{fig:p-dimer}
\end{center}
\end{figure}

Using \eqref{E-sol}, \eqref{S-sol} and \eqref{D-sol} we verify the normalization 
\begin{equation}
\label{norm-2}
\pi_0(t) + \pi_1(t) + \pi_2(t) = 1
\end{equation}
providing a consistency check. The plots of $\pi_0, \pi_1, \pi_2$ are shown in Fig.~\ref{fig:p-dimer}.

\subsection{Coverings by trimers ($\ell = 3$)}
\label{subsec:trimers}

The density $E_m$ now evolves according to
\begin{equation}
\label{Emt-3}
\dot E_m =-(m+2)E_m, \qquad m\geq 1
\end{equation}
from which
\begin{equation}
\label{Em-3-sol}
E_m(t)=e^{-(m+2)t}
\end{equation}
The density of voids is found from \eqref{V-E} and \eqref{Em-3-sol}
\begin{equation}
\label{Vm-3-sol}
V_m(t)=e^{-(m+2)t}(1-e^{-t})^2
\end{equation}

The density of empty sites is 
\begin{equation}
\label{empty-3-sol}
\pi_0(t) \equiv E_1(t) = e^{-3 t}
\end{equation}

The challenge is to determine the density $\pi_1$ of single-covered sites, the density $\pi_2$ of double-covered sites, and the density $\pi_3$ of triple-covered sites. It proves useful to consider 
\begin{equation}
\label{M-def:trimers}
\textsf{M} = \pi_2 + 2\pi_3
\end{equation}
This quantity grows according to the rate equation 
\begin{equation}
\label{M-eq}
\frac{d \textsf{M}}{dt} = 6 \sum_{m\geq 1} V_m
\end{equation}
which is derived similarly to \eqref{D-eq} for the dimer covering. Indeed, one notices that the change of $\textsf{M}$ occurs after a trimer lands on one or two covered sites at the boundary of a void. In the former case
\begin{equation}
\label{B1}
\bullet\circ\,\circ \longrightarrow  \blacktriangle\bullet\bullet
\end{equation}
where $\blacktriangle$ denotes a double-covered site. The void is assumed to contain at least two empty sites, that is, $m\geq 2$. In the process \eqref{B1}, the quantity $\textsf{M}$ undergoes the $0\to 1$ change. Combining with a similar deposition event on the right boundary of the void, we obtain a $2V_m$ contribution from voids of length $m\geq 2$.

If a new trimer covers two sites on the boundary of a void, the possible changes are
\begin{subequations}
\begin{align}
\label{B2-a}
&\bullet\bullet\,\circ \longrightarrow  \blacktriangle\blacktriangle\bullet\\
\label{B2-b}
&\blacktriangle\bullet\circ \longrightarrow  \clubsuit\blacktriangle\bullet
\end{align}
\end{subequations}
where $\clubsuit$ denotes a triple-covered site. The change of $\textsf{M}$ is $0\to 2$ in the process \eqref{B2-a} and  $1\to 3$ in the process \eqref{B2-b}.  Combining with a similar deposition event on the right boundary of the void, we obtain a $4V_m$ contribution from voids of length $m\geq 1$.

In the exceptional case of $m=1$, there is a single deposition event 
\begin{equation}
\label{V1-3}
\bullet\circ\,\bullet \longrightarrow  \blacktriangle\bullet \blacktriangle
\end{equation}
covering two boundaries of the void. The change of $\textsf{M}$ is $0\to 3$, yielding the gain term $2V_1$ similar to the gain term $2V_m$ coming from \eqref{B1} and its mirror image for voids of length $m\geq  2$. Therefore the overall rate of change of $\textsf{M}$ is indeed $2V_m+4V_m=6V_m$ summed over all $m\geq 1$ as stated in Eq.~\eqref{M-eq}. 

Inserting \eqref{Vm-3-sol} into Eq.~\eqref{M-eq} and integrating we obtain
\begin{subequations}
\begin{equation}
\label{M-3-sol}
\pi_2(t) + 2\pi_3(t)  =  \tfrac{1}{2} - 2 e^{-3t} + \tfrac{3}{2} e^{-4t} 
\end{equation}
Using normalization condition, $\pi_0+ \pi_1+\pi_2+ \pi_3 = 1$, together with \eqref{empty-3-sol} we arrive at
\begin{equation}
\label{SDT-3-sol}
 \pi_1(t)+\pi_2(t) + \pi_3(t) = 1-e^{-3t}
\end{equation}
\end{subequations}

\begin{figure}[ht]
\begin{center}
\includegraphics[width=0.456\textwidth]{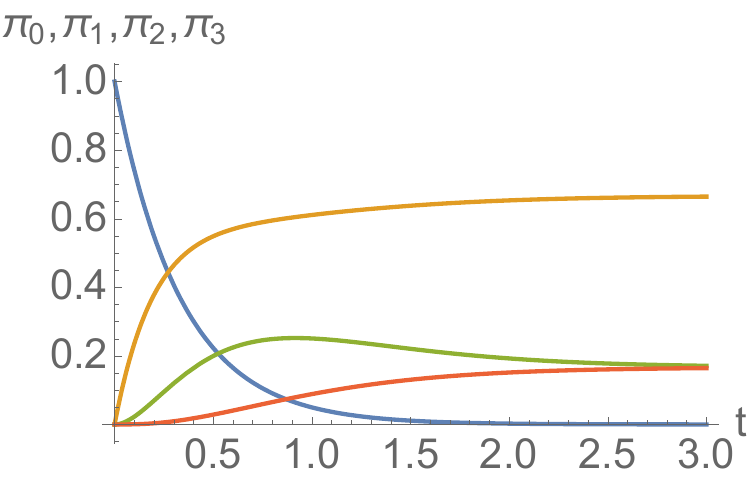}
\caption{The RSC by trimers. Shown is the evolution of the densities $\pi_0, \pi_1, \pi_2, \pi_3$ of uncovered, single-covered, double-covered, and triple-covered sites. The density $\pi_0(t)$ decays exponentially, Eq.~\eqref{empty-3-sol}. Other densities are the combinations of three exponents. The conjectural step \eqref{p3} has led to expressions \eqref{p123} plotted in the figure. These densities satisfy $\pi_1(t)>\pi_2(t)>\pi_3(t)$ for all $t>0$. }
\label{fig:p0123-3}
  \end{center}
\end{figure}

To fix the densities $\pi_1(t), \pi_2(t), \pi_3(t)$ one needs another relation supplementing \eqref{M-3-sol}--\eqref{SDT-3-sol}. So far we haven't derived such a relation. Instead, we rely on a less solid theoretical approach.  An analogy with the dimer case suggests that the densities $\pi_1(t), \pi_2(t), \pi_3(t)$ are the linear combinations of the exponents $\{e^{-2t}, ~e^{-3t}, ~e^{-4t}\}$. In Appendix~\ref{app:small} we show that it suffices to make a guess
\begin{equation}
\label{p3}
\pi_3 = a_0 + a_2 e^{-2t}  + a_3  e^{-3t} + a_4  e^{-4t}
\end{equation}
and then one can fix all parameters in \eqref{p3}. This leads [see Appendix~\ref{app:small} for details] to the following (conjectural) expressions for the densities (see also Fig.~\ref{fig:p0123-3}):
\begin{subequations}
\label{p123}
\begin{align}
\label{p1-3}
& \pi_1 =  \tfrac{2}{3} -  e^{-2t}  + \tfrac{7}{3}  e^{-3t} - 2 e^{-4t}\\
\label{p2-3}
& \pi_2 =  \tfrac{1}{6}  + 2e^{-2t}   - \tfrac{14}{3} e^{-3t} + \tfrac{5}{2} e^{-4t}\\
\label{p3-3}
& \pi_3 =  \tfrac{1}{6}  - e^{-2t}   + \tfrac{4}{3} e^{-3t} - \tfrac{1}{2} e^{-4t}
\end{align}
\end{subequations}

In the $t\to\infty$ limit, the densities \eqref{p123} approach to the final densities given by \eqref{coverage-trimers}. Equations  \eqref{p1-3} and \eqref{p3-3} predict that the densities $\pi_1$ and $\pi_3$ are increasing functions of time, while \eqref{p2-3} gives the density $\pi_2$ with maximum at $t_*=\ln(5/2)=0.91629\ldots$.

\subsection{Coverings by $\ell-$mers with arbitrary $\ell$}
\label{subsec:ell-mers}

In the general case of arbitrary $\ell$, the density $E_m$ obeys
\begin{equation}
\label{Emt-gen}
\dot E_m =-(m+\ell-1)E_m, \qquad m\geq 1
\end{equation}
from which
\begin{equation}
\label{Em-gen-sol}
E_m(t)=e^{-(m+\ell-1)t}
\end{equation}
Therefore the density of empty sites is 
\begin{equation}
\label{empty-gen-sol}
\pi_0(t) \equiv E_1(t) = e^{-\ell t}
\end{equation}
and  the density of voids is
\begin{equation}
\label{Vm-gen-sol}
V_m(t)=e^{-(m+\ell -1)t}(1-e^{-t})^2
\end{equation}

The quantity
\begin{equation}
\label{M-def}
\textsf{M} = \sum_{k=2}^\ell (k-1)\pi_k
\end{equation}
counting multiply covered sites with proper multiplicity generalizes $\textsf{M}=\pi_2$ for dimers and $\textsf{M}=\pi_2+2\pi_3$ for trimers. In the general case 
\begin{equation}
\label{M-eq-gen}
\frac{d \textsf{M}}{dt} = \ell(\ell-1) \sum_{m\geq 1} V_m
\end{equation}
extending Eqs.~\eqref{D-eq} and \eqref{M-eq} to arbitrary $\ell$. By inserting \eqref{Vm-gen-sol} into \eqref{M-eq-gen} we obtain 
\begin{equation*}
\dot {\textsf{M}} = \ell(\ell-1) (1-e^{-t})  e^{-\ell t}
\end{equation*}
from which
\begin{subequations}
\begin{equation}
\label{M-gen-sol}
\textsf{M}(t)  =  \tfrac{\ell-1}{\ell+1} - (\ell-1) e^{-\ell t} + \tfrac{\ell(\ell-1)}{\ell+1} e^{-(\ell +1)t} 
\end{equation}
The normalization condition gives
\begin{equation}
\label{norm-ell}
\sum_{k=1}^\ell \pi_k=1-e^{-\ell t}
\end{equation}
\end{subequations}

Equations \eqref{M-gen-sol}--\eqref{norm-ell} give two relations for $\pi_1,\ldots,\pi_\ell$. One needs $\ell-2$ extra relations to fix $\pi_k(t)$.

\section{Model B}
\label{sec:mod-B}

In this section, we consider the RSC by $\ell$-mers using model B dynamics accepting only those deposition events where overlap with previously deposited $\ell-$mers is at most $\lfloor \frac{\ell}{2}\rfloor$. For any odd $\ell$, this rule implies that the middle site of the incoming $\ell$-mer lands into the uncovered site. 

In the case of the dimer coverage, models A and B are identical. Models A and B differ starting from the trimer coverage. We derive complete results when $\ell=3,4,5$ and then discuss the general situation. 

\subsection{Coverings by trimers ($\ell = 3$)}
\label{subsec:trimers-B}

The densities of voids satisfy
\begin{equation}
\label{Vm-B} 
\frac{d V_m}{dt}  = -m V_m + 2 \sum_{n\geq m+2} V_n 
\end{equation}
The form of \eqref{Vm-B} suggests to seek the solution in the form
\begin{equation}
\label{Vm-ansatz}
V_m(t) = \Phi(t)\,e^{-mt}
\end{equation}
Substituting this ansatz into \eqref{Vm-B} we obtain
\begin{equation}
\label{Phi}
\frac{d \Phi}{dt} = 2\Phi\,\frac{e^{-2t}}{1-e^{-t}}
\end{equation}
Integrating this equation subject to the initial condition
\begin{equation}
\label{IC}
\lim_{t\to 0} t^{-2}\Phi(t) =1
\end{equation}
ensuring that the lattice was initially empty we fix $\Phi(t)$ and determine the void distribution
\begin{equation}
\label{Vm-3-B}
V_m(t) = (1-e^{-t})^2\,e^{2e^{-t}-2}\,e^{-mt}
\end{equation}

\begin{figure}[ht]
\begin{center}
\includegraphics[width=0.456\textwidth]{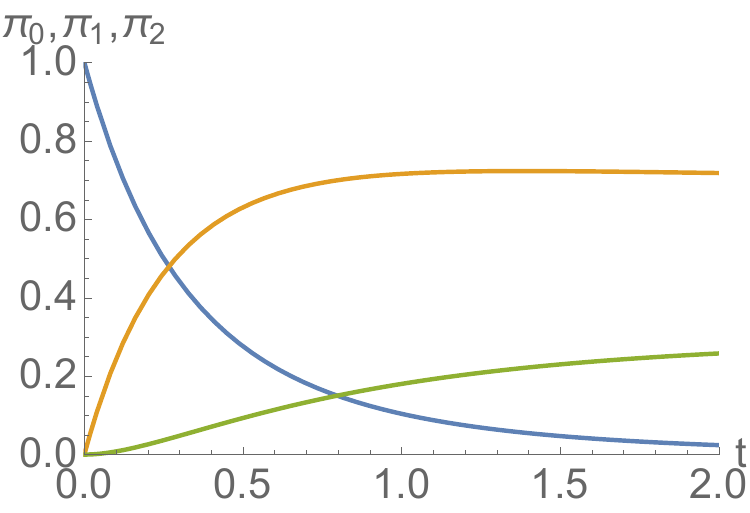}
\caption{The RSC by trimers, model B. The densities $\pi_0, \pi_1, \pi_2$ given by \eqref{p0-B}, \eqref{p1-B}, and \eqref{p2-B} are plotted. The fraction $\pi_1$ of single-covered sites always exceeds the fraction $\pi_2$ of double-covered sites: $\pi_1(t)>\pi_2(t)$ for all $t>0$. The fraction of double-covered sites always grows, while the fraction of single-covered sites reaches maximum, $\pi_1(t_*)\approx0.72313$, at $t_*\approx 1.386294$. }
\label{fig:p012-3}
\end{center}
\end{figure}

The fraction of uncovered sites is
\begin{equation}
\label{p0-B}
\pi_0(t) = \sum_{m\geq 1}mV_m= e^{-t}\,e^{2e^{-t}-2}
\end{equation}
The fraction of double-covered evolves according to the rate equation 
\begin{equation}
\label{p2-B-eq}
\frac{d \pi_2}{dt}  = 2 \sum_{m\geq 1} V_m = 2(1-e^{-t})\,e^{-t}\,e^{2e^{-t}-2}
\end{equation}
Integrating \eqref{p2-B-eq} subject to the initial condition $\pi_2(0)=0$ we obtain 
\begin{equation}
\label{p2-B}
\pi_2(t)  =\frac{1-(3-2e^{-t})\,e^{2e^{-t}-2}}{2}
\end{equation}
The fraction of single-covered sites $\pi_1=1-\pi_0-\pi_2$ is
\begin{equation}
\label{p1-B}
\pi_1(t)  =\frac{1+(3-4e^{-t})\,e^{2e^{-t}-2}}{2}
\end{equation}
This fraction reaches maximum at $t_*\approx 1.38629$. The plots of $\pi_0, \pi_1, \pi_2$ are shown in Fig.~\ref{fig:p012-3}. 

\subsection{Coverings by $4-$mers}
\label{subsec:4-B}

The densities of voids satisfy
\begin{equation}
\label{Vm-B-4}
\frac{d V_m}{dt}  = -(m+1) V_m + 2 \sum_{n\geq m+2} V_n 
\end{equation}
The form of \eqref{Vm-B-4} suggests to seek the solution in the form
\begin{equation}
\label{Vm-ansatz-4}
V_m(t) = \Phi(t)\,e^{-(m+1)t}
\end{equation}
One arrives at the same equation \eqref{Phi} for $\Phi(t)$ as before, and hence the void distribution is 
\begin{equation}
\label{Vm-4-B}
V_m(t) = (1-e^{-t})^2\,e^{2e^{-t}-2}\,e^{-(m+1)t}
\end{equation}

The fraction of uncovered sites is
\begin{equation}
\label{p0-B-4}
\pi_0(t) =\sum_{m\geq 1}mV_m=e^{-2t}\,e^{2e^{-t}-2}
\end{equation}
The fraction of double-covered evolves according to the rate equation 
\begin{equation}
\label{p2-B-eq-4}
\frac{d \pi_2}{dt}  = 6 \sum_{m\geq 1} V_m = 6(1-e^{-t})\,e^{-2t}\,e^{2e^{-t}-2}
\end{equation}
from which
\begin{equation}
\label{p2-B-4}
\pi_2(t)  =3(1-e^{-t})^2\,e^{2e^{-t}-2}
\end{equation}
leading to the announced final fractions \eqref{cov-4-B}. 

\begin{figure}[ht]
\begin{center}
\includegraphics[width=0.456\textwidth]{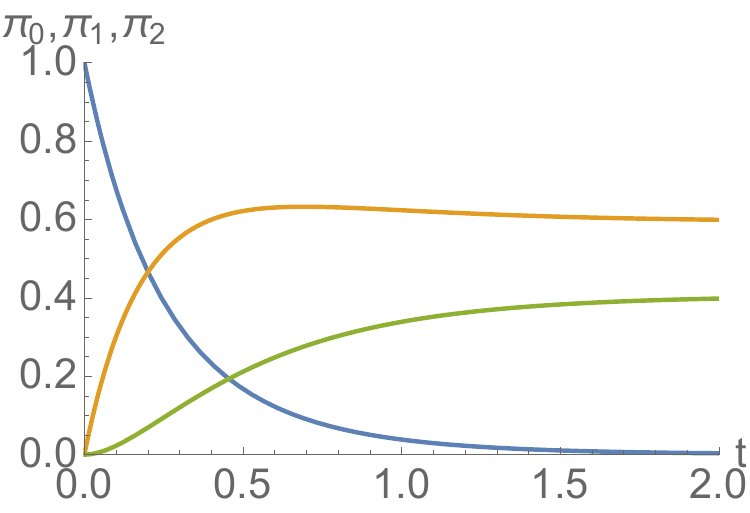}
\caption{The RSC by $4-$mers, model B. Shown is the evolution of the densities $\pi_0, \pi_1, \pi_2$ of uncovered, single-covered, and double-covered sites. The density $\pi_0(t)$ is given by Eq.~\eqref{p0-B-4}, the density $\pi_2(t)$ is given by Eq.~\eqref{p2-B-4}. }
\label{fig:p012-4}
\end{center}
\end{figure}

The plots of $\pi_0, \pi_1, \pi_2$ are shown in Fig.~\ref{fig:p012-4}. 
The fraction of single-covered sites reaches maximum $\pi_1(t_*)=1-1/e$ at $t_*=\ln 2$. This maximum is more visible than in the case of the covering by trimers (cf. Fig.~\ref{fig:p012-3} and Fig.~\ref{fig:p012-4}).

\subsection{Coverings by $5-$mers}
\label{subsec:5-B}

The densities of voids satisfy
\begin{equation}
\label{Vm-B-5}
\frac{d V_m}{dt}  = -m V_m + 2 \sum_{n\geq m+3} V_n 
\end{equation}
Seeking the solution in the form \eqref{Vm-ansatz} we simplify the infinite system \eqref{Vm-B-5} to a single differential equation
\begin{equation}
\frac{d \Phi}{dt} = 2\Phi\,\frac{e^{-3t}}{1-e^{-t}}
\end{equation}
Solving this equation we  determine the void distribution
\begin{equation}
\label{Vm-5-B}
V_m(t) = (1-e^{-t})^2\,e^{e^{-2t}+2e^{-t}-3}\,e^{-mt}
\end{equation}

\begin{figure}[ht]
\begin{center}
\includegraphics[width=0.456\textwidth]{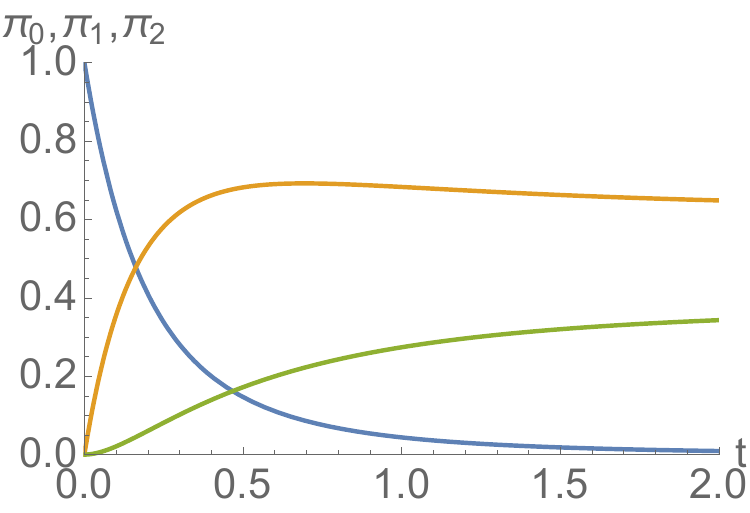}
\caption{The RSC by $5-$mers, model B. Shown is the evolution  of the densities $\pi_0, \pi_1, \pi_2$ of uncovered, single-covered, and double-covered sites. The density $\pi_0(t)$ is given by Eq.~\eqref{p0-B-5}, the density $\pi_2(t)$ is given by Eq.~\eqref{p2-B-5}.}
\label{fig:p012-5}
\end{center}
\end{figure}

The fraction of uncovered sites is
\begin{equation}
\label{p0-B-5}
\pi_0(t)=\sum_{m\geq 1}mV_m=e^{-t}\,e^{e^{-2t}+2e^{-t}-3}
\end{equation}
The fraction of double-covered sites varies according to the rate equation 
\begin{eqnarray}
\label{p2-B-5-eq}
\frac{d \pi_2}{dt}  &=& 4V_1+6 \sum_{m\geq 2} V_m \nonumber\\
&=& 2e^{-t}\,(1-e^{-t})\,(2+e^{-t})\,e^{e^{-2t}+2e^{-t}-3}
\end{eqnarray}
Integrating \eqref{p2-B-5-eq} subject to $\pi_2(0)=0$ we obtain 
\begin{eqnarray}
\label{p2-B-5}
\pi_2(t) &=& -1+e^{-t}\,e^{e^{-2t}+2e^{-t}-3}\nonumber\\
&+&\frac{5\sqrt{\pi}}{2\,e^4}\,\{\text{Erfi}[2]-\text{Erfi}[1+e^{-t}]\}
\end{eqnarray}
leading to the announced final fraction \eqref{cov-5-B}. The fraction $\pi_1$  of single-covered sites reaches maximum at the same time $t_*=\ln 2$ as for the RSC by $4-$mers. The magnitudes of the maxima are different, and generally $\pi_0, \pi_1, \pi_2$ differ (cf. Fig.~\ref{fig:p012-4} and Fig.~\ref{fig:p012-5}). 

\subsection{General case}
\label{subsec:gen-B}

The general case of arbitrary $\ell$ can be similarly treated. For instance, if $\ell$ is even, $\ell = 2p$ where $p$ is an arbitrary positive integer, the void distribution satisfies
\begin{equation}
\label{Vm-B-p}
\frac{d V_m}{dt}  = -(m+1) V_m + 2 \sum_{n\geq m+p} V_n 
\end{equation}
The solution has the form  \eqref{Vm-ansatz-4} with 
\begin{subequations}
\begin{equation}
\label{Phi-Pi}
\Phi(t) = (1-e^{-t})^2\,e^{\Pi_p(t)}
\end{equation}
where we shortly write
\begin{equation}
\label{Pi}
\Pi_p(t) = 2\sum_{q=1}^{p-1} (-1)^q \binom{p-1}{q}\frac{ (1-e^{-t})^q}{q}
\end{equation}
\end{subequations}

The fraction of uncovered sites is
\begin{equation}
\label{p0-gen-B}
\pi_0 = e^{-2t}\,e^{\Pi_p(t)}
\end{equation}
The fraction of double-covered sites varies according to the rate equation 
\begin{equation}
\label{p2-B-eq-p}
\frac{d \pi_2}{dt} = \sum_{k=1}^{p-2}(k+1)(2p-k)V_k +  p(p+1)\sum_{k\geq p-1}V_k
\end{equation}

Recalling that $V_k(t)=\Phi(t)\,e^{-(k+1)t}$, and using $\Phi$ defined by \eqref{Phi-Pi}--\eqref{Pi}, one can compute the sums  in the right-hand side of Eq.~\eqref{p2-B-eq-p} and then integrate to find $\pi_2(t)$. Thus we obtain a solution in the general case when $\ell = 2p$ is an arbitrary positive even integer. The solution is rather formal as the integral cannot be expressed through known special functions when $\ell\geq 8$.

\section{RSC of the line}
\label{sec:line}

Packing the line with randomly placed sticks is a classical ``car parking" problem going back to Renyi \cite{Renyi58}. Covering of an interval or a ring by sticks is also an old subject \cite{Stevens39,Domb47,Dvoretzky56,Flatto,Mandelbrot72,Shepp72a,Shepp72b,Siegel78,Siegel79,Holst,Shepp87}. In Secs.~\ref{sec:Dimers-Fluct}--\ref{sec:L-mers}, we studied the discrete models, viz. the RSCs of intervals by $\ell-$mers, and showed how to probe the full counting statistics of congested configurations. Extending such an analysis to the case of continuous random covering is left to the future. Here we consider the dynamics of the covering of an infinite line.  In the lattice case, there are two natural versions of the RSC process, model A and model B (see Secs.~\ref{sec:inf-time}--\ref{sec:mod-B}). In the continuous case, these two models are defined as follows:
\begin{itemize}
\item If a deposition attempt leads to an increase in coverage, it is accepted (model A).
\item Accepted deposition events are such that the center of incoming stick lands in the previously uncovered region (model B).
\end{itemize}
The evolution behaviors arising in continuous versions of models A and B are very different as we show below. 

\subsection{Model A}

Suppose that sticks are deposited on the line according to the rules of model A. Without loss of generality we set the length of sticks and the deposition rate to unity. The density $E(x,t)$ of empty intervals of length $x$ evolves according to
\begin{equation}
\label{Ext}
\frac{\partial E(x,t)}{\partial t}=-(x+1)E(x,t)
\end{equation}
from which
\begin{equation}
\label{Ex-sol}
E(x,t)=e^{-(x+1)t}
\end{equation}
The uncovered fraction of the line shrinks according to 
\begin{equation}
\label{p0-A}
\pi_0(t) \equiv E(0,t) = e^{-t}
\end{equation}

The density of the voids of length $x$ is found from the  continuous analog of \eqref{V-E}, viz.
\begin{equation}
\label{V-E-x}
V(x,t) = \frac{\partial^2 E(x,t)}{\partial x^2} 
\end{equation}
yielding
\begin{equation}
\label{Vx-sol}
V(x,t)=t^2\,e^{-(x+1)t}
\end{equation}

The fractions $\pi_k(t)$ of the line covered $k$ times are non-trivial for all $k\geq 0$. We already know $\pi_0=e^{-t}$, and the sum rule
\begin{equation}
\label{norm:p}
\sum_{k\geq 0} \pi_k(t) = 1
\end{equation}
implied by normalization. The initial condition is
\begin{equation}
\pi_k(0) = \delta_{k,0}
\end{equation}

As in the lattice case [cf. Eq.~\eqref{M-def}] it proves useful to consider the quantity
\begin{equation}
\label{M-def-cont}
\textsf{M}(t) = \sum_{k\geq  2}(k-1)\pi_k(t)
\end{equation}
accounting for all covered fractions with appropriate multiplicity. The choice of the multiplicity implies that in each deposition event  the quantity \eqref{M-def-cont} increases with rate equal to the overlap with already covered part of the line. We arrive at the evolution equation for $\textsf{M}(t)$
\begin{equation}
\label{M-cont}
\frac{d {\textsf{M}}}{dt} = \int_0^\infty dx\, V(x,t)=2t\, e^{-t}
\end{equation}
from which 
\begin{eqnarray}
\label{M-sol-line}
\textsf{M}(t)=2 -2(1+t)e^{-t} 
\end{eqnarray}

The challenge is to determine $\pi_k(t)$ with $k\geq 1$.  The limiting (jammed) values $\pi_k(\infty)$ are particularly interesting. The above results imply two constraints
\begin{equation}
\label{Sum-12}
\sum_{k\geq  1}\pi_k(\infty) = 1, \qquad \sum_{k\geq  1}(k-1)\pi_k(\infty) = 2
\end{equation}

The simplest guess is that $\pi_k(\infty)$  is a geometric distribution: $\pi_k(\infty)=A q^k$. If true, the sum rules \eqref{Sum-12} fix the distribution to 
$\pi_k(\infty) = 2^{k-1}/3^k$.

\subsection{Model B}

The density $V(x,t)$ of the voids of length $x$ varies in time according to  
\begin{equation}
\label{Vxt} 
\frac{\partial V(x,t)}{\partial t}=-x V(x,t)+2\int_{x+\frac{1}{2}}^\infty dy\,V(y,t)
\end{equation}
The uncovered fraction is
\begin{equation}
\label{p0-cont}
\pi_0(t)=\int_0^\infty dx\,x V(x,t)
\end{equation}

The line is initially empty. Therefore
\begin{subequations}
\begin{equation}
\label{IC-1} 
V(x,t=0)=0
\end{equation}
for all $x>0$. Using \eqref{p0-cont} we additionally deduce  
\begin{equation}
\label{IC-2}
\lim_{t\to 0}\int_0^\infty dx\,x V(x,t)=1
\end{equation}
\end{subequations}
supplementing \eqref{IC-1}.

The form of Eq.~\eqref{Vxt} suggests an exponential dependence of $V(x,t)$ on $x$, viz. $V(x,t)\propto e^{-xt}$. Thus we seek the solution in the form 
\begin{equation}
\label{exp}
V(x,t)=e^{-xt}\,t^2\Phi^2(t)
\end{equation}
(Writing the time-dependent factor as $t^2\Phi^2(t)$ is convenient for the analysis.) Substituting \eqref{exp} into the governing equation \eqref{Vxt} we obtain 
\begin{equation}
\label{Phi-eq}
\dot \Phi=-\Phi\,t^{-1}\big(1-e^{-t/2}\big)
\end{equation}
The initial condition \eqref{IC-1} is manifestly obeyed by the ansatz \eqref{exp}. The supplementary initial condition \eqref{IC-2} yields $\Phi(0)=1$. Solving \eqref{Phi-eq} subject to $\Phi(0)=1$ we find $\Phi(t)$. The void density \eqref{exp} becomes
\begin{equation}
\label{Vxt-sol}
V(x,t)=e^{-xt}\,t^2 \mathcal{E}(t)
\end{equation}
where 
\begin{equation}
\label{E-def}
\mathcal{E}(t) \equiv \exp\left[-2\int_0^{2t} du\,\frac{1-e^{-u}}{u}\right]
\end{equation}
The function $\mathcal{E}(t)$ appears in many formulas, e.g., it gives the uncovered fraction
\begin{equation}
\label{p0-line}
\pi_0(t)=\mathcal{E}(t)
\end{equation}
Using the asymptotic
\begin{eqnarray*}
\int_0^z du\,\frac{1-e^{-u}}{u}=\ln z +\gamma + z^{-1}e^{-z}+\ldots, 
\end{eqnarray*}
where $\gamma=0.577215\ldots$ is the Euler constant, we find that the uncovered fraction of $\mathbb{R}$ decays  as
\begin{equation}
\label{p0-asymp}
\pi_0(t)\simeq C\,t^{-2} \qquad\text{with} \quad
C=\frac{1}{4e^{2\gamma}}
\end{equation}
when $t\gg 1$.

\begin{figure}[ht]
\begin{center}
\includegraphics[width=0.456\textwidth]{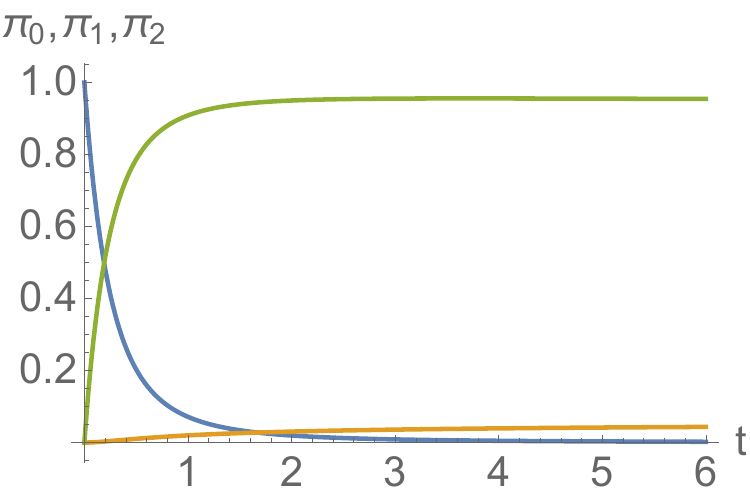}
\caption{Covering the line according to the rules of model B. Shown is the evolution of the densities $\pi_0, \pi_1, \pi_2$. The densities $\pi_1$ and $\pi_2$ given by Eqs.~\eqref{p1-line}--\eqref{p2-line} are increasing functions of time satisfying $\pi_1(t)>\pi_2(t)$ for all $t>0$. The density $\pi_0(t)$ decays according to Eq.~\eqref{p0-line}. }
\label{fig:p12-cont}
\end{center}
\end{figure}

In the realm of model B, the triple coverage is impossible in one dimension: $\pi_k=0$ for $k\geq 3$. We know $\pi_0$, Eq.~\eqref{p0-line}. To compute $\pi_2$, we notice that the double-covered fraction increases according to 
\begin{eqnarray}
\label{p2-eq}
\frac{d \pi_2}{dt} &=& \int_0^\frac{1}{2} dx\, x(1-x)V(x,t)\nonumber \\
&+& \frac{1}{4}\int_\frac{1}{2}^\infty dx\, V(x,t)          
\end{eqnarray}
Plugging \eqref{Vxt-sol}  into \eqref{p2-eq} and integrating we determine $\pi_2$; normalization $\pi_0+\pi_1+\pi_2=1$ then fixes $\pi_1$:
\begin{subequations}
\begin{align}
\label{p1-line}
\pi_1& = 1- \mathcal{E}(t) - \int_0^t d\tau\,\mathcal{E}(\tau)\left[1-2\,\frac{1-e^{-\tau/2}}{\tau}\right]\\
\label{p2-line}
\pi_2& = \int_0^t d\tau\,\mathcal{E}(\tau)\left[1-2\,\frac{1-e^{-\tau/2}}{\tau}\right]
\end{align}
\end{subequations}
see Fig.~\ref{fig:p12-cont}. The final densities are
\begin{equation}
\label{p12-B}
\pi_1(\infty) \approx 0.94553, \qquad \pi_2(\infty) \approx 0.054469
\end{equation}
The final densities \eqref{p12-B} show that model B leads to random covering close to optimal as the double-covered fraction constitutes less than $6\%$. 
The approach to these final densities is algebraic:
\begin{equation}
\label{p12-asymp}
\pi_1(t)-\pi_1(\infty)  \simeq \pi_2(\infty)-\pi_2(t) \simeq C\,t^{-1}
\end{equation}
with the same amplitude $C$ as in \eqref{p0-asymp}.

\section{Coverings the space}
\label{sec:space}

Here we analyze a random covering of $\mathbb{R}^d$ by balls of unit radius. We begin with model B exhibiting more interesting behavior than model A. According to the definition of model B, in a successful deposition attempt, the center of an incoming ball lands in the previously uncovered region. We now demonstrate that the uncovered volume fraction $\pi_0(t)$ vanishes algebraically in the long-time limit. 

To derive the asymptotic behavior of $\pi_0(t)$ we follow the same reasoning as in the case of random sequential adsorption (RSA). Namely, it is intuitive \cite{Feder80,Pomeau80,Swendsen81} to expect that in the large time limit, the uncovered space is essentially a collection of small uncovered patches, holes in short. Holes are small compared to balls. In two dimensions, for instance, the sides of each hole are almost straight (as they are circular arcs tiny compared to the radii of the balls). It seems intuitively obvious that almost all holes are triangles; the probability to have a hole with more than three sides becomes negligible in the large time limit. Similarly, holes resemble tetrahedra when $d=3$, and $d-$simplexes in the general case. These details are less relevant than two more intuitively obvious properties. The first is that we can roughly characterize each hole by its linear size $\ell$. The second is that holes are effectively independent in the long time limit. 

Denote by $c(\ell, t)$ the density of holes per linear size. The volume of a hole scales as $\ell^d$, and the total deposition rate is proportional to volume. Thus the density of well-separated and effectively non-interacting (in the large time limit) holes satisfies $dc/dt\sim -\ell^d c$ from which
\begin{equation}
\label{clt}
c(\ell,t)\sim \exp(-\ell^d t)
\end{equation}
The uncovered volume fraction is therefore
\begin{equation}
\label{vol-B}
\pi_0(t)\sim \int d\ell\,\,\ell^d c(\ell,t)\sim \int d\ell\,\,\ell^d e^{-\ell^d t} \sim t^{-1-1/d}
\end{equation}
This derivation is a straightforward generalization of an argument \cite{Feder80,Pomeau80,Swendsen81} yielding the rate of approach to the jammed state in RSA. 

In one dimension, \eqref{vol-B} agrees with the exact asymptotic \eqref{p0-asymp}. The amplitude predicted by \eqref{p0-asymp} is beyond the heuristic approach that has led to \eqref{vol-B}. Two other non-trivial fractions $\pi_1$ and $\pi_2$ are given by \eqref{p1-line} and \eqref{p2-line} for model B in one dimension. For model B in $d$ dimensions, the non-trivial fractions are $\pi_k$ with $0\leq k\leq z_d$, i.e., $\pi_k=0$ for $k>z_d$. Thus $z_1=2$ in one dimension. The threshold values $z_d$ are unknown. Conjecturally, $z_d$ are related to the kissing numbers $K_d$, defined as a maximal number of non-overlapping unit spheres that can touch a common unit sphere \cite{Rogers,GL87,Conway}. Namely, $z_d=K_d$ apart from special dimensions $d=2,8,24$ and perhaps a few more where $z_d=K_d-1$. Even if the relation between $z_d$ and $K_d$ is correct,  it is not very useful as the kissing numbers are known \cite{VDW,Leech,Levenshtein79,Odlyzko79,Ziegler,Musin08} only when $d=1,2,3,4,8,24$. 

For model A, the fractions $\pi_k$ are positive for all $k$. We haven't computed $\pi_k$ even in one dimension; the only exception is $\pi_0(t)$, see \eqref{p0-A}. For model A in arbitrary dimension, the asymptotic decay of $\pi_0(t)$ can be computed more accurately than for model B. The heuristic picture of the long-time evolution remains the same, but the region where the center of an incoming ball can land approaches the ball of unit radius. The asymptotic form of the rate equation for $c(\ell,t)$ is $dc/dt\simeq -V_d c$ where we set the deposition rate to unity and $V_d= \pi^{d/2}/\Gamma(1+d/2)$ is the volume of the ball of unit radius. The uncovered volume fraction decays exponentially 
\begin{equation}
\label{vol-A}
\pi_0(t)\sim e^{-V_d t}
\end{equation}
in contrast to the algebraic decay \eqref{vol-B} for model B.

\section{Concluding Remarks}
\label{sec:CR}

Our definition of the RSC is an analog of the RSA. Among other random coverings, we mention the procedure introduced by Matheron \cite{Matheron68,Matheron},  known as visible confetti or the dead leaves model (DLM), see \cite{Serra,Kendall99,DLM-Mumford,DLM-Roueff,DLM-Galerne,DLM-Penrose,Klenke22}. In this model, leaves fall at random onto the ground long enough, so the ground is completely covered. The visible parts of leaves on the ground tessellate $\mathbb{R}^d$. The two-dimensional  DLM has received considerable attention in applications, e.g., in image modeling and material science \cite{DLM-Mumford,DLM-Roueff,DLM-Galerne}. In contrast to models A and B where the deposition depends on the successful previous  events, the deposition rule in the DLM is random. However, the focus in the DLM is on the visible part of leaves that is non-trivial already in one dimension. It would be interesting to apply the methods we used in the analysis of models A and B to the DLM. 

An interesting class of models concerns complete coverings of lattices. The most famous is the dimer covering problem solved for planar lattices in Refs.~\cite{Kasteleyn61,Fisher61}. The complete coverings of the triangular lattice by triangular trimers are also understood \cite{Nienhuis99,Nienhuis01}. These models are very different in spirit, namely, they belong to equilibrium statistical physics. 

Intriguing questions about densest packings and least dense coverings by balls concern the $d\to\infty$ behavior.  The upper and lower bounds for the fraction of the covered space in the packing problem and the average coverage in the covering problem are old \cite{Rogers,Conway,TS10,Cohn16,Parisi20}. The bounds are exponentially separated, e.g., for the fraction of the covered space the lower and upper bounds scale as $2^{-d}$ and $2^{-0.5990 d}$, respectively. Some problems are tractable analytically in the $d\to\infty$ limit  \cite{Frisch85,Frisch87,Kurchan16,Charbonneau17}. Finding the asymptotic average coverage of $\mathbb{R}^d$ by the RSC process (version B) may improve an upper bound for the average coverage. For the RSA, one would get an improved lower bound for the fraction of the covered space.

\appendix
\section{Model A: Covering by trimers}
\label{app:small}

Here we consider the RSC of the one-dimensional lattice by trimers using the rules of model A, i.e., accepting the deposition of a trimer if at least one previously uncovered site gets covered. The RSC by dimers is completely solved, namely the fractions $\pi_0, \pi_1, \pi_2$ are known. The RSC of trimers 
is the first lattice system that is not completely solved. Several exact results have been derived Sec.~\ref{subsec:trimers}, yet some basic quantities defy analytical treatment. For instance, we have determined the density $\pi_0$ of uncovered sites, Eq.~\eqref{E-sol}, while the densities $\pi_1$, $\pi_2$ and $\pi_3$ of single-covered, double-covered and triple-covered sites remain unknown. Here we derive expressions \eqref{p1-3}--\eqref{p3-3} for these densities relying on a rather natural conjecture.

Expressing $\pi_1$ and $\pi_2$ through $\pi_3$ ensuring that \eqref{M-3-sol} and \eqref{SDT-3-sol} are obeyed we obtain
\begin{subequations}
\begin{align}
\label{p13}
&\pi_1(t) =  \tfrac{1}{2} + e^{-3t} - \tfrac{3}{2} e^{-4t} + \pi_3(t)\\
\label{p23}
&\pi_2(t) =  \tfrac{1}{2} - 2 e^{-3t} + \tfrac{3}{2} e^{-4t} -2 \pi_3(t)
\end{align}
\end{subequations}
Hence it suffices to find $\pi_3$. 

In the dimer case, the density $\pi_2$ of the most occupied sites is a linear combination of $\{1, e^{-2t}, e^{-3t}\}$, cf. \eqref{D-sol}. In the trimer case, $\pi_3$ accounts for the most occupied sites and we {\em guess} that it is a linear combination of $\{1, e^{-2t}, e^{-3t}, e^{-4t}\}$; the necessity of the additional exponent is evident from \eqref{M-3-sol}. This leads to the conjectural form \eqref{p3} with four parameters. Expanding \eqref{p3} and matching with the small time behavior
\begin{equation}
\label{p3-small}
\pi_3 = \tfrac{2}{3}t^3+O(t^4)
\end{equation}
fixes the parameters in \eqref{p3} and gives \eqref{p3-3}. Inserting \eqref{p3-3} into \eqref{p13}--\eqref{p23} gives \eqref{p1-3}--\eqref{p2-3}.

To derive the asymptotic behavior \eqref{p3-small} we first notice that the rules of the RSC imply 
\begin{equation}
\label{pk-small}
\pi_k(t) = A_k(\ell) t^k + O\big(t^{k+1}\big)
\end{equation}
Thus Eq.~\eqref{p3-small} asserts that $A_3(3)=\frac{2}{3}$. To deduce this result we first show how to derive two other amplitudes, $A_1(3)=3$ and $A_2(3)=3$. 

In the small time limit,  almost all deposited objects are isolated trimers. Their density is $t$, and therefore $\pi_1\simeq 3t$. This asymptotic can be confirmed by expanding \eqref{p13} and taking into account that $\pi_3= O(t^3)$. 

The double-covered sites are generated via the deposition of trimers overlapping with already present isolated trimers:
\begin{subequations}
\begin{align}
\label{2-a}
&\cdots\circ\bullet\bullet\bullet\,\circ\circ\circ \cdots \longrightarrow  \cdots\circ\bullet\bullet\blacktriangle\bullet\bullet\circ \cdots\\
\label{2-b}
&\cdots\circ\bullet\bullet\bullet\,\circ\circ\circ \cdots \longrightarrow  \cdots\circ\bullet\,\blacktriangle\,\blacktriangle\bullet\circ\circ \cdots
\end{align}
\end{subequations}
The asymptotic contribution of the process \eqref{2-a} [resp. \eqref{2-b}] to production of $\pi_2$ is $t$   [resp. $2t$].  Accounting  for the symmetric events on the other boundary of the isolated trimer gives $\frac{d\pi_2}{dt}\simeq 6t$, so $\pi_2\simeq 3t^2$. This asymptotic can be confirmed by expanding \eqref{p23} and recalling that $\pi_3=O(t^3)$. 

The triple-covered sites are produced (in the leading order) by the process 
\begin{equation}
\label{p3:gain}
\cdots\circ\bullet\blacktriangle\blacktriangle\bullet\circ\circ \cdots \longrightarrow \cdots\circ\bullet\blacktriangle\clubsuit\blacktriangle\bullet\circ \cdots
\end{equation}
The density of the pattern on the left in \eqref{p3:gain} is $t^2$. This follows from our previous analysis: The patterns on the right in \eqref{2-a} and \eqref{2-b} arise with asymptotically equal densities which should be equal to $t^2$ since  $\pi_2\simeq 3t^2$. Combining the contribution of the process \eqref{p3:gain} and its mirror version yields $\frac{d\pi_3}{dt}\simeq 2t^2$ leading to Eq.~\eqref{p3-small}. 

\bigskip\noindent
I am grateful to Henk Hilhorst for helpful remarks and to a referee who informed me about the dead leaves model.

\bibliography{references-packing}

\end{document}